\documentclass[11pt]{article}

\usepackage{graphicx}
\usepackage{url}
\usepackage{amsmath}
\usepackage{amsfonts}
\usepackage{verbatim}
\usepackage{amssymb}
\usepackage{amsthm}
\usepackage{color}
\usepackage{float}
\usepackage[english]{babel}
\usepackage[autostyle]{csquotes}
\usepackage{mathtools}
\usepackage[usenames,dvipsnames]{xcolor}
\usepackage{mdwlist}
\usepackage{caption}
\usepackage{subcaption}
\usepackage{setspace}
\usepackage{bbm}
\usepackage{epstopdf}
\usepackage{cite}
\usepackage{pdflscape}
\usepackage{wrapfig}
\usepackage{empheq}
\usepackage{hyperref}
\hypersetup{
    colorlinks=true,
    linkcolor=blue,
    filecolor=magenta,      
    urlcolor=blue,
    citecolor=blue,
}

\newtheorem{theorem}{Theorem}

\newtheorem{lemma}{Lemma}

\usepackage{amsthm}
\usepackage[top=1in, bottom=1in, left=1in, right=1in]{geometry}

\setcaptionmargin{0.25in}

\usepackage{hyperref}

\title{Determining the source of period-doubling instabilities in spiral waves}
\author{Stephanie Dodson\thanks{Division of Applied Mathematics, Brown University, Providence, RI: \texttt{stephanie\_dodson@brown.edu}}
\and Bj\"{o}rn Sandstede\thanks{Division of Applied Mathematics, Brown University, Providence, RI: \texttt{bjorn\_sandstede@brown.edu}}}
\date{\today}

\begin{document}


\maketitle

\begin{abstract}
\noindent Spiral wave patterns observed in models of cardiac arrhythmias and chemical oscillations develop alternans and stationary line defects, which can both be thought of as period-doubling instabilities. These instabilities are observed on bounded domains, and may be caused by the spiral core, far-field asymptotics, or boundary conditions. Here, we introduce a methodology to disentangle the impacts of each region on the instabilities by analyzing spectral properties of spiral waves and boundary sinks on bounded domains with appropriate boundary conditions. We apply our techniques to spirals formed in reaction-diffusion systems to investigate how and why alternans and line defects develop. Our results indicate that the mechanisms driving these instabilities are quite different; alternans are driven by the spiral core, whereas line defects appear from boundary effects. Moreover, we find that the shape of the alternans eigenfunction is due to the interaction of a point eigenvalue with curves of continuous spectra.
\end{abstract}

\textbf{Keywords:} spiral waves, reaction-diffusion systems, stability, alternans, period-doubling

\textbf{AMS subject classifications:} 35B36, 35K57


\section{Introduction}

Systems of oscillatory and excitable media frequently express spiral wave patterns. Spiral waves are observed in laboratory settings in chemical oscillations in the Belousov-Zhabotinsky reaction \cite{Zaikin:1970te,Winfree:1972ud} and cell signaling in slime molds \cite{Robertson:1981uy}, and have been associated with arrhythmic heart rhythms \cite{Wiener:1946,Winfree:1994cz,PerezMunuzuri:1991gt}. These systems support rigidly rotating spirals with constant shape, but transitions to complex dynamics and unstable spirals are common.
 
\begin{figure}[H]
\centering
        \includegraphics[width=0.8\textwidth]{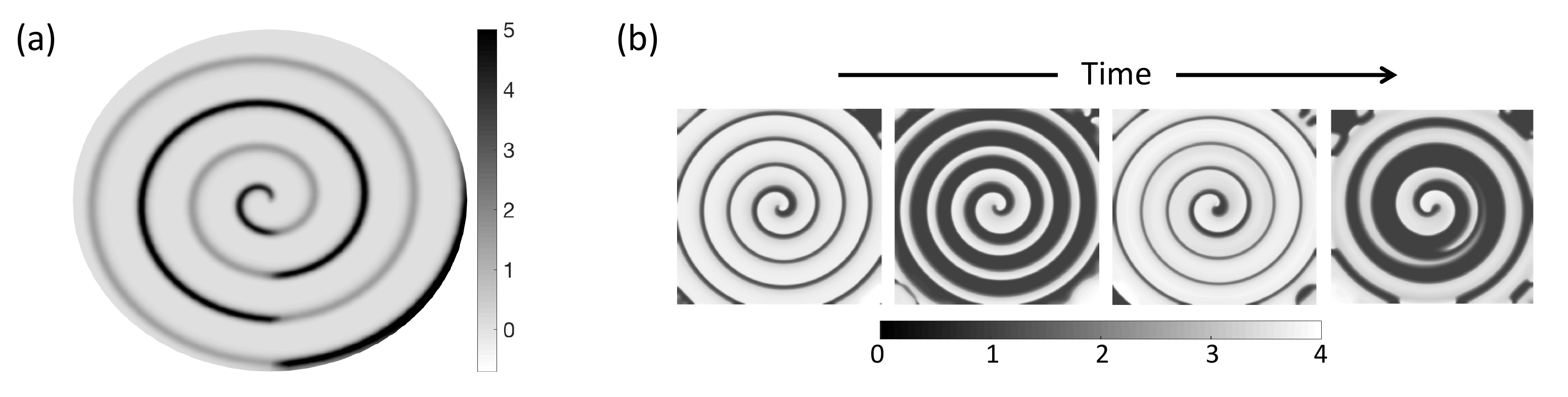}
    \caption{(a) Stationary line defect in the $w$-component of the R\"{o}ssler system. (b) Time evolution of alternans instability in the $u$-component of the Karma model on a square of side length 16cm with homogeneous Neumann boundary conditions. System parameters as defined in Section~\ref{sec:models}.  }\label{fig:alternans_line_defect}
\end{figure}

In cardiac dynamics, accelerated tachycardiac rhythms have been linked to electrical activity organized as rotating spiral wave patterns on the surface of the heart. The transition from tachycardiac to fibrillation is believed to be initiated by spiral wave breakup \cite{Rosenbaum:1994bt,Pastore:1999dq}. Clinical studies indicate a primary driver of breakup is conduction block following a long-short temporal modulation of the action potential duration, in what is known as the alternans instability. Alternans are visible on electrocardiograms and have become a clinical warning sign of sudden cardiac death \cite{Rosenbaum:1994bt,Pastore:1999dq}. In spiral waves, alternans physically corresponds to variation in spiral band width (Figure~\ref{fig:alternans_line_defect}). For a detailed review of spiral waves in cardiac dynamics, see the review article \cite{Alonso:2016et} and recent results published in the special issue \cite{Cherry:2017je}.

Spirals are also produced and studied in chemical oscillations, for example the Belousov-Zhabotinsky reaction \cite{Zaikin:1970te,Winfree:1972ud}. In these systems, spirals are experimentally observed to form stationary line defects \cite{Yoneyama:1995cn,Park:2006hz} (Figure~\ref{fig:alternans_line_defect}a), which have been reproduced in numerical simulations \cite{Goryachev:1998ky,Sandstede:2007jw}. Across the defect lines, wave amplitudes are out of phase (Figure~\ref{fig:alternans_line_defect}). Both alternans and line defects lead to a spiral wave with twice the period of the original planar wave.
 
Nonlinear reaction-diffusion systems qualitatively capture transitions to complex meandering, drifting, and the period-doubled line defects and alternans patterns. In these systems, planar spiral waves are stationary solutions in a rotating polar coordinate frame and converge to one-dimensional periodic travelling waves away from the core. Stability and bifurcations can be studied by considering the spectra of the operator obtained by linearizing the nonlinear system about the spiral wave solution. The spectrum consists of isolated eigenvalues and a set determined by the operator in the far-field limit. 

Bounded domains are of interest in applications to cardiac dynamics and laboratory experiments. Neumann boundary conditions naturally represent lower conductance tissue separating regions of the heart or the physical walls of containers. Mathematically, on finite domains planar spiral waves are truncated and matched with a boundary sink, which adds extra structure to the spiral wave and alters the spectrum of the linear operator \cite{Fiedler:2000,Sandstede:2000ug,Sandstede:2002ht}. The boundary sink itself directly contributes an additional set of eigenvalues, and the finite domain modifies the spectrum associated with the far-field dynamics. Furthermore, radial growth in the eigenfunctions is permitted, and those that would not be integrable on the full plane now emerge as true eigenfunctions on the bounded domain \cite{Fiedler:2000,Sandstede:2000ug,Sandstede:2002ht}. These eigenvalues are associated with intrinsic properties of the spiral wave and are attributed to the spiral core. The spectrum of the operator on a bounded domain is therefore a union of three disjoint sets that are associated with instabilities from the far-field, boundary conditions, and core. Knowing which set unstable eigenvalues belong to provides information about how instabilities will manifest themselves on unbounded or bounded domains. 
 
Meander and drift instabilities are the result of a Hopf bifurcation originating from the core: the emerging dynamics is understood through actions of the symmetry group of translations and rotations on the plane and a center manifold reduction \cite{Barkley:1994vk,Sandstede:1999wl}. However, previous studies investigating alternans and line defects provide inconsistent and incomplete evidence for which spectral set the unstable eigenvalues belong to. 
 
Due to the clinical significance, the alternans instability has been a recent area of focus in the cardiac dynamics community. In single cells, alternans are widely attributed to a period-doubling instability observed in simple 1D maps \cite{Guevara:1984vu}. However, this condition, known as the restitution hypothesis, has received contradictory evidence \cite{Cytrynbaum:2002ca,Cherry:2004caa} and does not appear to be relevant for excitable tissues that support traveling waves. The formation and stability of alternans in waves propagating on a ring and line have been analyzed with kinematic descriptions \cite{Echebarria:2007ir} and through linear stability analyses \cite{Bar:1999ez,Echebarria:2002ca,Bar:2004jc,Bauer:2007dd}. Stability analysis in 1D predicts that alternans are the result of a Hopf bifurcation \cite{Gottwald:2008hi}, yet analysis of the 1D traveling waves cannot fully capture 2D features. 
 
Linear stability analysis of spirals on bounded domains in the Karma and Fenton-Karma models found a variety of unstable eigenmodes responsible for the formation of alternans \cite{Marcotte:2015ex, Marcotte:2016hn, Allexandre:2004cd}. In \cite{Marcotte:2015ex} and \cite{Marcotte:2016hn}, Marcotte and Grigoriev find that formation of alternans depends on the domain size. Furthermore, they determine that the alternans eigenmodes are not spatially localized near the core.  
 
The  rigorous analysis of spiral waves in \cite{Sandstede:2007jw} indicates period-doublings are initiated by a series of Hopf bifurcations with imaginary parts of the eigenvalues sitting robustly at multiples of half the spiral frequency. These Hopf eigenvalues may be induced by period-doubling of the far-field dynamics or boundary sinks. Stationary line defects are hypothesized to stem from bifurcations of the boundary sink, but no direct evidence supporting this claim was found in \cite{Sandstede:2007jw} .
 
The goal of this paper is to further investigate how and why these period-doubling like instabilities arise on bounded domains. Specifically, we seek to answer which spectral set the unstable eigenvalues belong to and gain a better understanding of how the spirals destabilize and on what domains the instabilities are relevant. To tackle this problem, we introduce a methodology for disentangling the contributions of each region by forming related patterns on domains whose spectra will contain eigenvalues arising from a subset of the resulting spectra. Three cases are considered and compared with essential and absolute spectra from wave trains: (1) a spiral on a bounded disk with Neumann boundary conditions to provide the full spectrum, (2) a boundary sink to demonstrate effects of boundary conditions, and (3) a spiral on a disk with non-reflecting boundary conditions to remove boundary eigenvalues.

We apply this methodology to reaction-diffusion models and discover that the mechanisms driving alternans and line defects are rather different. We find that line defects arise from the boundary sink and thus will only appear under the correct conditions on bounded domains. In contrast, alternans originate from instabilities associated with the spiral core and will develop independent of the domain. Furthermore, the spectral computations reveal that the structure of the alternans instability develops due to an interaction of an unstable point eigenvalue and curves of continuous spectra. Our results have important consequences for reproducing patterns such as line defects and provide justification for extending analysis of alternans from simple bounded disks to the complex geometry of the heart.

In the sections that follow, we begin by describing the mathematical set-up and notation of the reaction-diffusion models. We include a  review of relevant spectral properties for operators on the plane and how these properties are modified by bounded domains. Procedures used to compute the spirals and spectra are described in the methods section. Finally, we present the results of the analysis applied to the R\"{o}ssler and Karma models to study line defects and alternans, respectively.

\section{Models}  \label{sec:models}

Reaction-diffusion systems display a rich set of patterns and are commonly used to model systems in biology and nature. General planar reaction-diffusion systems are of the form
\begin{align} \label{eqn:gen_rxn_diff}
U_t = D \Delta U + F(U), \ \ \ U \in \mathbb{R}^n,  \ \  D \in \mathbb{R}^{n \times n},   \ \  x \in \mathbb{R}^2, 
\end{align}
where $U = (u_1, \dots, u_n)^T$ is a vector of species that diffuse at rates given by the nonnegative elements $\delta_i$ of the diagonal matrix $D$, and $\Delta$ is the Laplace operator. Kinetic reactions of the different species are captured by the typically nonlinear function $F(U)$.
 
Cardiac models range in complexity from biophysically detailed ion-channel models to simplified systems which capture qualitative features, with both categories falling under the reaction-diffusion framework. The Karma system is a two-variable reduction of the Noble ion-channel model \cite{Noble:1962tk} and was developed to be a simplified model that reproduces alternans \cite{Karma:1993tq,Karma:1994gb}. The model is given by
\begin{align} \label{eqn:karma}
u_t &= 1.1 \Delta u + 400 \left( -u + \left(1.5414 - v^4 \right) \left(1 - \tanh(u - 3)\right) \frac{u^2}{2} \right)\\
v_t &= 0.1 \Delta v + 4 \left( \frac{1}{1 - e^{-\mu_K}} \vartheta \left(u - 1 \right) - v \right), \nonumber
\end{align}
where the fast variable $u$ represents voltage and $v$ acts as a slower gating variable. As in \cite{Marcotte:2015ex, Marcotte:2016hn,Allexandre:2004cd}, we use the function $\vartheta(u) = \left(1 + \tanh(4 u) \right)/2$. Alternans are observed in this system when the real bifurcation parameter $\mu_K$ is increased above one \cite{Karma:1994gb}.
 
The R\"{o}ssler model is commonly used to study chaotic turbulence in chemical oscillations, and is known to produce spirals with line defects. This three-variable system is also of the general reaction-diffusion form and is given by
\begin{align} \label{eqn:rossler}
u_t &= 0.4 \Delta u - v - w\\
v_t & = 0.4 \Delta v + u + 0.2 v \nonumber \\
w_t & = 0.4 \Delta w + uw - \mu_R w + 0.2. \nonumber
\end{align}
Bifurcations to line defects are observed for parameter values $\mu_R > 3$ \cite{Goryachev:1998ky,Sandstede:2007jw}. 
 
Here, we write both models in a general form and define the bifurcation parameters to be $\mu_K$ and $\mu_R$, respectively. We remark that in (\ref{eqn:karma}) and (\ref{eqn:rossler}) we selected values for several parameters that are often allowed to vary. We refer to Table~\ref{table:params} in the appendix for the general form of these models.

\section{Review of Spiral Waves and their Spectral Properties}

In our context, periodic traveling waves, also referred to as wave trains, serve as building blocks of spiral waves. Therefore, we first consider the existence and stability properties of wave trains on $\mathbb{R}$ and their restriction to bounded domains before describing spiral waves on the plane and bounded disks. Further details can be found in \cite{Fiedler:2000,Sandstede:2000ug,Sandstede:2002ht,Kapitula:2013}.

\subsection{Wave trains and boundary sinks}

On $\mathbb{R}$, the reaction-diffusion system (\ref{eqn:gen_rxn_diff}) reduces to
\begin{align} \label{eqn:rxn_diff_1D}
U_t = D U_{xx} + F(U), \ \ x \in \mathbb{R}.
\end{align}
Wave trains are solutions to (\ref{eqn:rxn_diff_1D}) of the form $U(x,t) = U_{\infty}(\kappa x - \omega t)$ where $U_{\infty}$ is $2\pi$-periodic in its argument, so that $\kappa$ is the spatial wave number, and $\omega$ is the temporal frequency. In the traveling coordinate $\xi = \kappa x - \omega t$,  wave trains are stationary solutions of
\begin{align} \label{eqn:wt}
U_t = \kappa^2 D U_{\xi \xi} + \omega  U_{\xi}  + F(U), \ \ \ \ \xi \in \mathbb{R}.
\end{align}
Generically, wave trains arise as one-parameter families for which $\omega$ and $\kappa$ are connected by the nonlinear dispersion relation $\omega = \omega_*(\kappa)$ and the profile $U_{\infty}(\xi;\kappa)$ depends smoothly on $\kappa$. The group velocity of the wave train is defined from the nonlinear dispersion relation as $\text{c}_{\text{g}} = \frac{d \omega}{d \kappa}$: it is equal to the speed with which perturbations are transported along the wave train in the original laboratory frame. 
 
\begin{figure}[H]
\centering
        \includegraphics[width=0.7\textwidth]{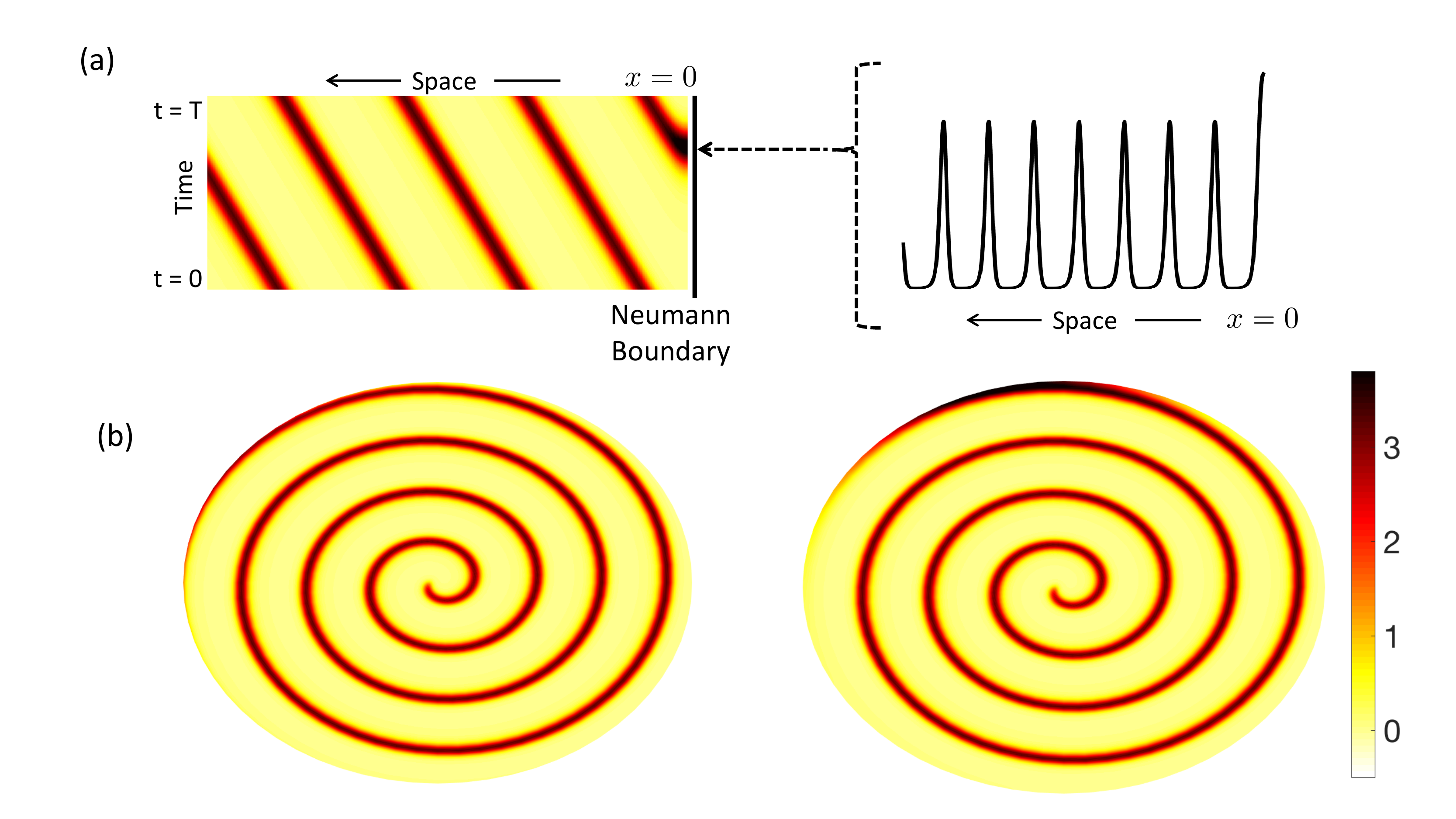}
    \caption{(a) Illustration of a boundary sink. The effect of Neumann boundary condition is highlighted by the temporal cross section on the right. (b) Comparison between planar spiral wave (left) and spiral wave on a bounded disk with Neumann boundary conditions (right). Both spirals shown on disks for comparison purposes. }\label{fig:boundary_sink_bounded_spiral}
\end{figure}

To prepare for our discussion of spiral waves on bounded domains, we introduce the concept of boundary sinks which connect wave trains of (\ref{eqn:rxn_diff_1D}) at $x = -\infty$ with a Neumann boundary condition at $x = 0$. We say that $U(x,t) = U_{\text{bdy}}(x,\omega t)$ where $U_{\text{bdy}}(x,\tau)$ is $2\pi$-periodic in $\tau$ and satisfies the following one-dimensional equation on the half line in the laboratory frame 
\begin{align} \label{eqn:bndry_sink}
& \omega U_{\tau} = D U_{xx}  + F(U), \ \ (x,t) \in (-\infty,0) \times S^1,\\
&U_x(0,\tau) = 0 , \ \ \tau \in  S^1 \nonumber 
\end{align}
and converges to a wave train $U_{\infty}(\kappa x - \omega t)$ with $\text{c}_\text{g} > 0$ as $x \rightarrow -\infty$ such that
\begin{align*}
\big| U_{\text{bdy}}(x,\cdot) - U_{\infty}(\kappa x - \cdot) \big|_{C^1(S^1)} \rightarrow 0.
\end{align*} 
 An example of a boundary sink is shown in Figure~\ref{fig:boundary_sink_bounded_spiral}. 

\subsection{Planar spiral waves and truncation to bounded disks}
We say that the reaction-diffusion system (\ref{eqn:gen_rxn_diff}) has a planar spiral wave solution of the form $U(x,t) = U_*(r,\phi - \omega t)$ where $(r,\phi)$ are polar coordinates if there exists an $\omega \in \mathbb{R}$ and a smooth function $\theta(r)$ with $\theta'(r) \rightarrow 0$ such that $U_*$ satisfies (\ref{eqn:gen_rxn_diff}) and 
\begin{align*}
|U_*(r,\cdot - \omega t) - U_{\infty}(\kappa r + \theta(r) + \cdot - \omega t) |_{C^1\left(S^1\right)} \rightarrow 0, \ \ \ \text{ as }  r \rightarrow \infty,
\end{align*}
where $U_{\infty}(\kappa r - \omega t)$ is a wave train with c$_{\text{g}} > 0$. Here, $\omega$ is the temporal rotational frequency and $\theta(r)$ acts as a phase correction to match solutions at the core with asymptotic wave trains. The spatial wave number $\kappa$ is selected by the spiral, and the wave train connects $\omega$ and $\kappa$ through the nonlinear dispersion relation $\omega = \omega_*(\kappa)$. Spiral waves are stationary solutions in the co-rotating polar frame $(r,\psi) = (r,\phi - \omega t)$
\begin{align} \label{eqn:gen_rxn_diff_polar}
U_t = D \Delta_{r,\psi} U + \omega  U_{\psi} + F(U).
\end{align}

Finite domains are physically and numerically realistic, and bounded disks in particular are common computational domains as they incorporate rotational symmetry properties of the spiral. When considered on $B_R(0)$, the disk of radius $R$ centered at the origin, planar spiral waves are truncated and solutions with positive group velocity emitted from the core are now matched with time $2\pi/\omega$-periodic boundary sinks $U_{\text{bdy}}(\xi,\tau)$. Spirals formed on $B_R(0)$ with homogeneous Neumann boundary conditions are stationary solutions of the system in the rotating polar frame
\begin{align} \label{eqn:bdd_spiral}
&U_t = D \Delta_{r,\psi} U + \omega_R U_{\psi} + F(U), \ \ (r,\psi) \in [0,R) \times S^1\\
&U_r(R,\psi) = 0, \ \ \psi \in S^1. \nonumber
\end{align}
The temporal frequency of the bounded spiral converges to that on the  infinite domain $\omega_R \rightarrow \omega_*$ as $R \rightarrow \infty$. An example of spiral waves on bounded disks and the entire plane is shown in Figure~\ref{fig:boundary_sink_bounded_spiral}. Note how the homogeneous Neumann boundary conditions influence the spiral near the outer boundary at the top of the spiral, similar to the boundary sink.

\subsection{Spectral stability of linear operators}

Here, we give a review of relevant spectral and stability concepts. First, general spectral definitions and terminology is defined. Then the spectra of the one-dimensional wave trains is described, followed by those of spirals on unbounded domains. In each case, we first consider the patterns on infinite domains and then describe how spectra are modified by bounded domains. 
 
The spectrum $\Sigma$ of a closed, densely defined linear operator $\mathcal{L}: X \rightarrow X$ on the Banach space $X$ is defined as 
\begin{align*}
\Sigma = \left\{ \lambda \in \mathbb{C} \ \big| \ (\mathcal{L} - \lambda): X \rightarrow X \text{ does not have a bounded inverse } \right\}.
\end{align*}
If $\mathcal{L}$ is Fredholm, the spectrum can be decomposed into the following disjoint sets \cite{Kapitula:2013}:
\begin{align*}
\Sigma = \Sigma_{\text{pt}} \ \cup \ \Sigma_{\text{Fred}} \ \cup \ \Sigma_{\text{FB}}
\end{align*}
where
\begin{align*}
&\Sigma_{\text{pt}} = \left \{ \lambda \in \mathbb{C} : \mathcal{L} - \lambda \text{ is Fredholm with index 0 but not invertible} \right\}\\
& \Sigma_{\text{FB}} = \left\{  \lambda \in \mathbb{C} : \mathcal{L} - \lambda \text{ is not Fredholm}  \right\}\\
&\Sigma_{\text{Fred}} =  \left\{ \lambda \in \mathbb{C} : \mathcal{L} - \lambda \text{ is Fredholm with non-zero index}  \right\}.
\end{align*}
The set $\Sigma_{\text{pt}}$ is referred to as the point spectrum and contains elements called eigenvalues, which are typically discrete.  Often, the essential spectrum $\Sigma_{\text{ess}}$ is defined by $\Sigma_{\text{ess}} = \Sigma / \Sigma_{\text{pt}} =  \Sigma_{\text{Fred}} \ \cup \ \Sigma_{\text{FB}}$. The contents of each set will depend on the operator and some of these sets may be empty. These sets are illustrated in Figure~\ref{fig:fredholm_spectral_sets}.

\subsection{Stability of wave trains}

Stability of wave trains in the co-moving coordinate frame $\xi = \kappa x - \omega t$ is analyzed by considering the spectrum of the operator 
\begin{align} \label{eqn:wt_lin_op}
\mathcal{L}^{\text{mv}}_{\infty}V  = \kappa^2 D  V_{\xi \xi} + \omega V_{\xi} + F_U(U_{\infty}) V
\end{align}
on $L^2(\mathbb{R})$ formed by linearizing (\ref{eqn:wt}) about the solution $U_{\infty}$. There are no non-trivial solutions $V \in L^2(\mathbb{R})$ in the kernel of $\mathcal{L}^{\text{mv}}_{\infty} - \lambda$ and the sets $\Sigma_{\text{pt}}$ and $\Sigma_{\text{Fred}}$ are empty, that is $\Sigma_{\text{pt}} = \Sigma_{\text{Fred}} = \emptyset$. The spectrum $\Sigma$ of $\mathcal{L}^{\text{mv}}_{\infty}$ consists only of Fredholm borders $\Sigma_{\text{FB}}$, which can be computed as follows. The linearization $F_U(U_{\infty})$ is $2\pi$-periodic, so by Floquet theory we seek non-trivial solutions to $\mathcal{L}^{\text{mv}}_{\infty} V = \lambda V$ of the form $V(\xi,t) =  e^{\nu \xi/\kappa} \bar{V}(\xi)$ with $\Bar{V}(\xi + 2\pi) = \bar{V}(\xi)$ for $\nu \in \mathbb{C}$ and obtain the relation
\begin{align} \label{eqn:wt_lin_disp_rel}
\mathcal{L}^{\text{mv}}_{\infty}(\lambda,\nu)\bar{V} =  D \left( \kappa \partial_{\xi} + \nu \right)^2 \bar{V} + \frac{\omega}{\kappa} \left( \kappa \partial_{\xi} + \nu \right) \bar{V}+ F_U(U_{\infty})\bar{V} - \lambda \bar{V} = 0
\end{align}
which connects the temporal eigenvalues $\lambda$ and spatial Floquet exponents $\nu$. Therefore, the spectrum $\Sigma$ of $\mathcal{L}^{\text{mv}}_{\infty}$ is given by 
\begin{align*}
\Sigma_{wt} := \left\{ \lambda \in \mathbb{C} : \exists \ \nu \in i \mathbb{R} \text{ and non-trivial } 2\pi \text{-periodic } \bar{V}(\xi) \text{ so that } \mathcal{L}^{\text{mv}}_{\infty}(\lambda, \nu)\bar{V} = 0 \ \forall \xi \in \mathbb{R} \right\}.
\end{align*}
It can be shown that $\Sigma_{wt}$ is the union of smooth curves of the form $\lambda = \lambda_{\infty}(\nu)$ with $\nu = i \gamma \in i \mathbb{R}$, which are often referred to as linear dispersion curves. For each fixed $\lambda \in \mathbb{C}$, equation~(\ref{eqn:wt_lin_disp_rel}) admits finitely many Floquet exponents $\nu \in \mathbb{C}$, and at least one Floquet exponent crosses the imaginary axis as $\lambda$ crosses through a spectral curve.
 
\begin{figure}[H]
\centering
 \includegraphics[width=0.8\textwidth]{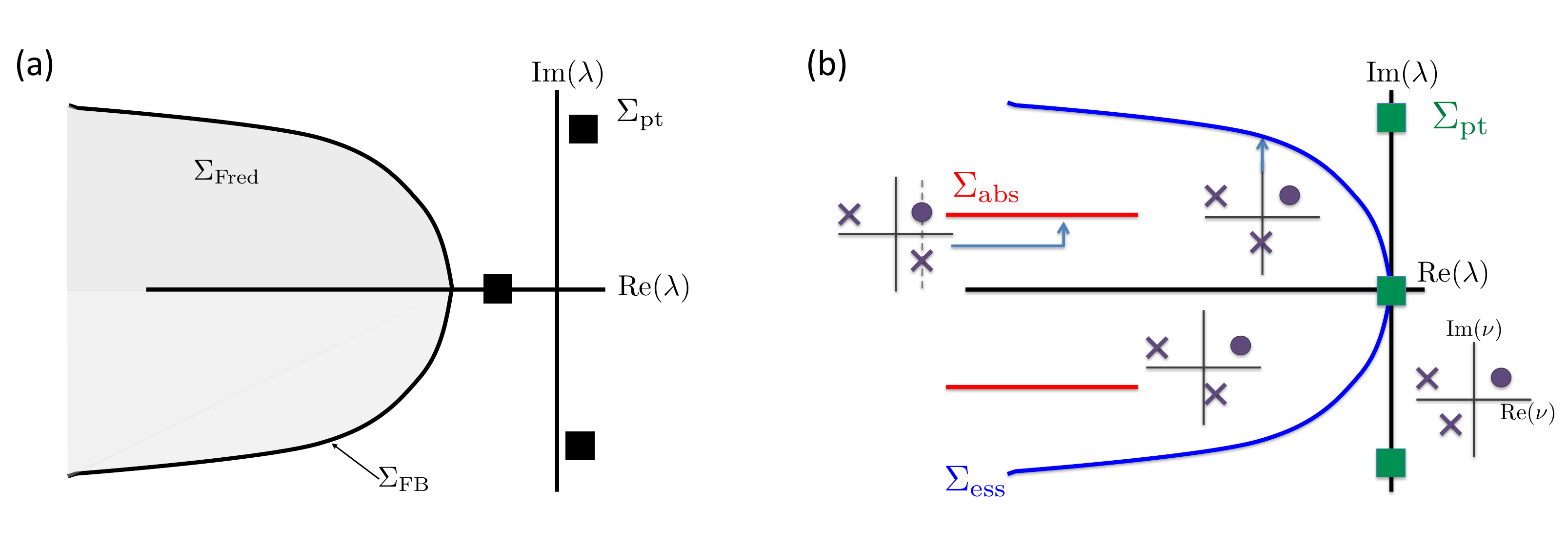}
\caption[Illustration of spectral sets]{ (a) Illustration of spectral sets for a general linear operator. In this example, the operator $\mathcal{L} - \lambda$ is Fredholm with index 0 in the unshaded region. Point spectrum (squares) can be found in this region and are defined for $\lambda$ that make the operator not invertible. Moving to the left, the operator $\mathcal{L}-\lambda$ becomes not Fredholm on the curve $\Sigma_{\text{FB}}$ and is then Fredholm with non-zero index in the shaded region.  (b) Illustration of spectra for spiral wave linear operator. Insets show distribution of spatial eigenvalues, with dots (crosses) indicating $\nu$ with initially Re$(\nu) >0$ ($<0$) for $\lambda \gg 1$. Crossing the essential spectrum from right to left results in one spatial eigenvalue crossing the imaginary axis and the real part becoming positive. On absolute spectra curves, the eigenvalue that crossed the imaginary axis aligns with one from the initial positive set. } \label{fig:fredholm_spectral_sets}
\end{figure}

In the laboratory frame, the linearized equation is
\begin{align} \label{eqn:rxn1d_labFrame}
V_t &= D V_{xx} + F_U \left(U_{\infty}(\kappa x - \omega t) \right)V.
\end{align}
Functions of the form $V(x,t) = e^{\lambda t} e^{\nu x} \bar{V}(\kappa x - \omega t)$ for non-trivial $2\pi$-periodic $\bar{V}(\kappa x - \omega t)= \bar{V}(\xi)$ satisfy equation~(\ref{eqn:rxn1d_labFrame}) if and only if
\begin{align} \label{eqn:disp_rel_lab}
\mathcal{L}_{\infty}^{\text{lab}}(\lambda,\nu) \bar{V} =   D \left( \kappa \partial_{\xi} + \nu \right)^2 \bar{V} + \omega \bar{V}_{\xi} + F_U \left(U_{\infty}(\xi) \right) \bar{V} - \lambda \bar{V} = 0,
\end{align}
which defines essential spectrum curves $\lambda = \lambda_{\text{lab}}(\nu)$ for $\nu \in i \mathbb{R}$. We note that the essential spectra in the co-moving (\ref{eqn:wt_lin_disp_rel}) and laboratory frames (\ref{eqn:disp_rel_lab}) are different: comparing $\mathcal{L}^{\text{mv}}_{\infty}(\lambda,\nu)$ to $\mathcal{L}^{\text{lab}}_{\infty}(\lambda,\nu)$, we see that the spectral curves are related via \cite{Sandstede:2000ut,Sandstede:2007jw}
\begin{align} \label{eqn:ess_spec_wt_spiral_relation}
\lambda_{\text{lab}}(\nu) = \lambda_{\infty}(\nu) - \frac{\omega}{\kappa} \nu + i \omega \ell, \ \ \ell \in \mathbb{Z}.
\end{align}
Essential spectra computed in the laboratory frame have vertical periodic branches parameterized by $\ell \in \mathbb{Z}$, which arise from Floquet ambiguity. Additionally, for each fixed $\lambda \in \mathbb{C}$ in (\ref{eqn:disp_rel_lab}), there are infinitely many $\nu \in \mathbb{C}$ and non-trivial $2\pi$-periodic functions $\bar{V}$ such that $\mathcal{L}^{\text{lab}}_{\infty}(\lambda,\nu)\bar{V} = 0$. We order the spatial eigenvalues $\nu_j$ for fixed $\lambda \gg 1$ by their real part,
\begin{align} \label{eqn:spatial_eigenvalues_ordered}
\cdots \leq  \text{Re}(\nu_{-j-1}) \leq \text{Re}(\nu_{-j}) \leq \cdots \leq  \text{Re}(\nu_{-1}) < 0 <   \text{Re}(\nu_{1}) \leq \cdots \leq \text{Re}(\nu_{j}) \leq \text{Re}(\nu_{j+1}) \leq \cdots  
\end{align}
so that $\text{Re}(\nu_{-1}) < 0 <   \text{Re}(\nu_{1})$. 
Upon crossing an essential spectrum curve (Fredholm border), at least one spatial eigenvalue crosses the imaginary axis. Figure~\ref{fig:fredholm_spectral_sets} shows curves of essential spectrum with insets indicating the distribution of the spatial Floquet eigenvalues $\nu$. 

\subsection{Stability of planar and bounded spiral waves}
Stability of planar spirals can be determined similarly to the one-dimensional case by considering the spectrum of the operator formed by linearizing~(\ref{eqn:gen_rxn_diff_polar}) about $U_*(r,\psi)$ 
\begin{align} \label{eqn:planar_lin_op_full}
\mathcal{L}_* V = D \left(\partial_{rr} + \frac{1}{r} \partial_r + \frac{1}{r^2} \partial_{\psi \psi} \right) V + \omega V_{\psi} + F_U(U_*(r,\psi)) V
\end{align}
and considering the operator $\mathcal{L}_*$ on $L^2(\mathbb{R}^2)$ with domain $H^2(\mathbb{R}^2)$. Here, the spectrum $\Sigma$ contains $\lambda = 0$ and $\lambda = \pm i \omega$ which arise from rotational and translational symmetries of the planar spiral.
 
Features of $\Sigma_{\text{ess}}$ depend purely on asymptotic properties of the spiral \cite{Sandstede:2000ug,Sandstede:2000ut}. In the formal limit $r \rightarrow \infty$, the linear operator $\mathcal{L}_*$ becomes
\begin{align} \label{eqn:planar_lin_op}
\tilde{\mathcal{L}}_{*} = D  \partial_{rr} + \omega \partial_{\psi} + F_U(U_{\infty}).
\end{align}
Eigenfunctions in the far-field limit take the form \cite{Sandstede:2001bb}
\begin{align} \label{eqn:efcn_form}
V(r,\psi) = e^{ \nu r} e^{i \ell \psi} \bar{V}(\kappa r + \psi) , \ \  \bar{V}(\xi + 2\pi) = \bar{V}(\xi),
\end{align}
where radial growth or decay is characterized by the real part of the spatial eigenvalue $\nu$ and $\bar{V}(\xi)$ is a periodic eigenfunction of the asymptotic wave train. Substitution into $\tilde{\mathcal{L}}_{*} V = \lambda V$ gives
\begin{align} \label{eqn:ff_spec}
\tilde{\mathcal{L}}_*(\lambda,\nu) \bar{V} =  D \left( \kappa \partial_{\xi} + \nu \right)^2 \bar{V} + \omega \bar{V}_{\xi} + F_U(U_{\infty}) \bar{V} - \lambda \bar{V}.
\end{align}
The Fredholm borders of the essential spectrum are defined by $\lambda = \lambda_*(\nu)$ for which one spatial eigenvalue is purely imaginary \cite{Sandstede:2000ug} and we see that the far-field spiral-wave operator reduces to the case of the laboratory frame wave train (\ref{eqn:disp_rel_lab}), that is $\lambda_*(\nu) = \lambda_{\text{lab}}(\nu)$. In general,
\begin{align*}
\Sigma_{\text{ess}} \left(\mathcal{L}_* \right)   = &\left\{  \lambda \in \mathbb{C} : \tilde{\mathcal{L}}_*(\lambda,\nu) \bar{V} = 0 \text{ has a non-trivial solution $\bar{V} \in H^2\left(S^1\right)$ for } \nu \in i \mathbb{R} \right\}\\
& \cup \left\{ \lambda \in \mathbb{C} : \mathcal{L} - \lambda_* \text{ is Fredholm with non-zero index } \right\} ,
\end{align*}
and this set is connected to the essential spectrum of periodic wave trains in the co-moving frame via relation (\ref{eqn:ess_spec_wt_spiral_relation}) \cite{Sandstede:2000ut}. Since $\nu\in i \mathbb{R}$, the mapping does not modify stability properties, and $\Sigma_{\text{ess}}\left(\mathcal{L}_{\infty}\right)$ and $\Sigma_{\text{ess}} \left(\mathcal{L}_* \right)$ destabilize under the same conditions. The additional vertical periodic branches at integer multiples of $i \omega$ are distinct branches for the spiral and no longer artifacts from Floquet theory as here far-field rotational symmetry implies that if $V(r,\psi)$ is an eigenfunction, then so is $e^{i\ell \psi}V(r,\psi)$. As in the laboratory frame, there exists infinitely many spatial eigenvalues $\nu$ for each temporal eigenvalues $\lambda$, which we will order by real part, as in (\ref{eqn:spatial_eigenvalues_ordered}).
%
 
The pertinent linear operator for spiral waves on the bounded disk $B_R(0)$ with Neumann boundary conditions is
\begin{align} \label{eqn:bdd_spiral_linOp}
&\mathcal{L}_{*,R}V = D \Delta_{r,\psi}V + \omega V_{\psi}  + F_U(U_{*,R}) V, \ \ (r,\psi) \in [0,R) \times S^1\\
&V_r(R,\psi) = 0, \ \ \ \psi \in S^1 \nonumber.
\end{align}
The spectrum $\Sigma_{*,R}$ of $\mathcal{L}_{*,R}$ contains only point spectrum as the operator $\mathcal{L}_{*,R} - \lambda$ is Fredholm with index zero for all $\lambda$. Naturally, we expect that the discrete eigenvalues in $\Sigma_{*,R}$ resemble the spectra of the planar spiral wave $\Sigma(\mathcal{L}_*)$, however, we will see that this is not true in general. Instead, eigenvalues in $\Sigma_{*,R}$ will converge to the union of three sets: the extended point spectrum $\Sigma_{\text{ext}}$, the absolute spectrum $\Sigma_{\text{abs}}$, and the spectrum of the boundary sink $\Sigma_{\text{bdy}}$. We describe these sets and their properties below.

 
Intuitively, the essential spectrum describes convective instabilities, in which growing perturbations are transported away to infinity \cite{Sandstede:2000ug}. When posed on a bounded disk, convective instabilities are no longer relevant. Instead, perturbations that grow in norm at every point in space become significant. These so-called absolute instabilities are captured in the limit $R\to\infty$ by the absolute spectrum, which is defined via the far-field linear dispersion relation $\tilde{\mathcal{L}}_*(\lambda,\nu)$ as
\begin{align*}
\Sigma_{abs} = \left \{ \lambda \in \mathbb{C} : \text{Re }\nu_{-1}(\lambda) = \text{Re }\nu_{1}(\lambda) \right \}.
\end{align*}
The absolute spectrum consists typically of curves that are parametrized by $\beta = |\text{Im }\nu_{-1} - \text{Im }\nu_1|$, where $\beta = 0$ at the end points. Elements of $\Sigma_{\text{abs}}$ do not correspond to eigenvalues of $\mathcal{L}_{*,R}$ but rather represent accumulation points of infinitely many discrete eigenvalues of $\mathcal{L}_{*,R}$ as the domain size $R$ goes to infinity \cite{Sandstede:2000ug}. Note that the absolute spectrum is still defined by the limiting operator $\tilde{\mathcal{L}}_{*}$ for the planar spiral wave, but it is generally distinct from the essential spectrum. The spatial eigenvalues corresponding to elements in the absolute spectrum are also illustrated in Figure~\ref{fig:fredholm_spectral_sets}.
 
To explain the extended point spectrum $\Sigma_{\text{ext}}$, we introduce spaces with exponential weight functions in polar coordinates on $\mathbb{R}^2$ with given weight $\eta \in \mathbb{R}$ in the radial direction via
\begin{align*}
L^2_{\eta}\left(\mathbb{R}^2 \right) := \{ u \in L^2_{\text{loc}} : |u|_{L^2_{\eta}} < \infty \}, \ \ \ |u|^2_{L^2_{\eta}} := \int_{\mathbb{R}^2} \big| u(x) e^{\eta |x|} \big|^2 dx.
\end{align*}
For every $\lambda \notin \Sigma_{\text{abs}}$, there exists an $\eta \in \mathbb{R}$ such that $\text{Re }\nu_{-1}(\lambda) < \eta < \text{Re }\nu_1(\lambda)$. Using exponentially weighted spaces, we can then define the extended point spectrum $\Sigma_{\text{ext}}$ via \cite{Sandstede:2000ug,Fiedler:2000}
\begin{align*}
\Sigma_{\text{ext}} = \big\{ \lambda \in \mathbb{C}\setminus\Sigma_{\text{abs}} : \mathcal{L}_* - \lambda \text{ is not boundedly invertible on } L^2_\eta \text{ where } \eta \\
\text{ is such that Re }\nu_{-1}(\lambda) < \eta < \text{Re }\nu_1(\lambda) \big\}.
\end{align*}
The  weight permits exponential radial growth of eigenfunctions up to rate $\eta$. Exponentially weighted norms are equivalent on bounded domains, therefore we can expect that elements in the extended point spectrum of planar spiral waves $\Sigma_{\text{ext}}(\mathcal{L}_*)$ persist as eigenvalues for the operator on each bounded disk, which we can attribute to instabilities caused by the core.

Finally, the boundary conditions may contribute additional point eigenvalues, which belong to the spectrum of the boundary sink $\Sigma_{\text{bdy}}$ defined above in equation~(\ref{eqn:bndry_sink}). The pertinent linearized operator is given by
\begin{align} \label{eqn:boundary_sink_plane}
&\mathcal{L}_{\text{bdy}} V = -\omega V_{\tau} + D V_{xx} + F_U(U_{\text{bdy}})V, \ \  (x,\tau) \in (-\infty,0) \times S^1\\
&V_x(0,\tau) = 0, \ \ \tau \in S^1. \nonumber
\end{align}
The relevant eigenvalues of the boundary sink are those that persist on finite domains \cite{Fiedler:2000,Sandstede:2000ug,Sandstede:2007jw}. Therefore, the extended point spectrum of the boundary sink provides the eigenvalues of the spiral caused by the boundary conditions, and we have $\Sigma_{\text{bdy}} = \Sigma_{\text{ext}} \left( \mathcal{L}_{\text{bdy}} \right)$.

We summarize our discussion in the following theorem from \cite[Theorem 2.5.5]{Fiedler:2000}, see \cite{Sandstede:2000ug,Sandstede:2019} for additional details. 
 
\begin{theorem}\label{thm:spectra_spirals} 
The spectrum $\Sigma\left( \mathcal{L}_{*,R} \right)$ of $\mathcal{L}_{*,R}$ on $L^2\left([0,R] \times [0,2\pi)\right)$ with Neumann boundary conditions converges: 
\begin{align*}
\Sigma \left(\mathcal{L}_{*,R} \right)\rightarrow \Sigma_{\text{abs}}\left( \mathcal{L}_* \right) \cup \Sigma_{\text{ext}}\left( \mathcal{L}_* \right) \cup \Sigma_{\text{ext}}\left(\mathcal{L}_{\text{bdy}} \right) \ \ \ \ \text{ as R } \rightarrow \infty.
\end{align*}
where $\mathcal{L}_*$ is the operator for a planar spiral wave on $L^2(\mathbb{R}^2)$ and $\mathcal{L}_{\text{bdy}}$ is the boundary sink operator. Convergence is uniform on bounded subsets of the complex plane in the symmetric Hausdorff distance. Moreover, the multiplicity of eigenvalues in the extended point spectrum is preserved; in contrast, the number of eigenvalues, counted with multiplicity, in any fixed open neighborhood of any point $\lambda \in \Sigma_{\text{abs}}\left( \mathcal{L}_* \right)$ converges to infinity as $R \rightarrow \infty$. 
\end{theorem}

Therefore, on bounded domains, eigenvalues fall into one of three sets: (1) extended point spectrum that persist under truncation, (2) eigenvalues converging to and emerging from the absolute spectrum, and (3) spectrum of the boundary sink.

We note that it was proved in \cite{Sandstede:2006ei} that, under certain conditions on the asymptotic equations, isolated eigenvalues in $\Sigma_{\text{ext}}$ may emerge from absolute spectrum branch points at predictable angles and destabilize prior to the absolute spectrum. The location of these isolated eigenvalues is predicted by including $1/r$ curvature terms into the asymptotic problem \cite{Sandstede:2006ei,Wheeler:2006eu}. 

\section{Methods}
Alternans and line defects are observed on bounded domains. To investigate whether unstable eigenvalues that generate these instabilities originate from $\Sigma_{\text{abs}}$, $\Sigma_{\text{ext}}$, or $\Sigma_{\text{bdy}}$, we consider spirals formed on three domains, each of which contains eigenvalues from a portion of the spectral sets.

The first domain is the standard bounded disk of radius $R$ with homogeneous Neumann boundary conditions, which we denote by $B_R(0)$. Here, spirals $U_{*,R}(r,\psi)$ are solutions of (\ref{eqn:bdd_spiral}), and the spectrum of the operator 
\begin{align}
\mathcal{L}_{*,R}V = D \Delta_{r,\psi} V + \omega V_{\psi} + F_U(U_*)V
\end{align}
on $L^2(B_R(0))$ with domain $\left\{V(x)\in H^2(B_R(0)): V_r(R,\cdot)=0 \right\}$ provides information about stability. From Theorem~\ref{thm:spectra_spirals}, the spectrum contains contributions from the core, the far field, and the boundary sink captured by the sets $\Sigma_{\text{ext}}$, $\Sigma_{\text{abs}}$, and $\Sigma_{\text{bdy}}$, respectively.
 
Core instabilities associated with $\Sigma_{\text{ext}}$ are analyzed by computing spirals on the bounded disk radius $R$ with non-reflecting boundary conditions. Non-reflecting boundary conditions mimic an infinite domain by matching the spiral to the asymptotic wave train on the boundary, allowing the spiral to naturally pass through it without interference. Non-reflecting spirals $U_{\text{nr}}(r,\psi)$ are solutions to
\begin{align} \label{eqn:asymp_spiral}
&0 = D \Delta_{r,\psi} U + \omega U_{\psi} + F(U), \ \  (r,\psi) \in [0,R) \times S^1 \nonumber \\
&0 =U_r - \kappa U_{\psi} , \ \ r = R, \ \psi \in  S^1,
\end{align}
where the boundary condition is obtained by taking derivatives of the asymptotic matching condition ${U_{*}(r,\psi)=U_{\infty}(\kappa r-\psi)}$. The linear operator for the non-reflecting spirals $\mathcal{L}_{\text{R,nr}}$ is
\begin{align}
\mathcal{L}_{\text{R,nr}}V = D \Delta_{r,\psi} V + \omega V_{\psi} + F_U(U_{\text{nr}})V
\end{align}
which acts on eigenfunctions in $ \left\{ V(x) \in H^2 \left( B_R(0)\right) : V_r(R,\psi) = \kappa V_{\psi}(R,\psi) \right\}$.
 
Finally, effects of the boundary and eigenvalues associated with $\Sigma_{\text{bdy}}$ are captured by direct computation of the boundary sinks. The time $2\pi/\omega$-periodic pattern $U_{\text{bdy}}$ is posed on a two-dimensional spatiotemporal domain $\Omega_{\text{bdy}} = [-L,0]  \times S^1$, with Neumann boundary conditions at $x = 0$ and $2\pi$-periodic boundary conditions in $\tau$. $U_{\text{bdy}}(x,\tau)$ is a solution to 
\begin{align}
&U_t{\tau}= \omega \left[D  U_{xx} + F(U) \right], \ \ (x,\tau) \in [-L,0) \times S^1\\
& U_x(0,\tau) = 0, \ \tau \in S^1.
\end{align}
 Stability of the boundary sink is given by considering the operator
\begin{align}
\mathcal{L}_{\text{bdy,L}} V = -\omega V_{\tau} + D V_{xx} + F_U(U_{\text{bdy}})V
\end{align}
on the space $\left\{ V(x,\tau) \in H^1([-L,0] \times S^1) : V_x(0,\cdot) = 0 \right\}$. Boundary sinks for the R\"{o}ssler and Karma models are shown in Figures~\ref{fig:fig4}a and~\ref{fig:fig6}b and will be discussed in further detail below.

Boundary sinks have far-field dynamics and boundary conditions, but lack the core: conversely, non-reflecting spirals contain the core but lack outer boundary effects. On each domain, stability properties are given by spectra of the operator linearized around the solution. Comparing the spectra of these three operators will indicate which region and spectral set is responsible for observed instabilities. We expect all operators to have eigenvalues aligning along the absolute spectrum due to the far-field dynamics, but the spectrum of $\mathcal{L}_{\text{bdy,L}}$ will not contain isolated core eigenvalues from $\Sigma_{\text{ext}}\left( \mathcal{L}_* \right)$, and $\mathcal{L}_{\text{R,nr}}$ will not have eigenvalues from the boundary sink. We remark that for non-reflecting boundary conditions in $\mathcal{L}_{\text{R,nr}}$ discrete eigenvalues from the far-field will still converge to the absolute spectrum. These expectations are summarized in the following lemma \cite{Sandstede:2019}.
\begin{lemma}
The spectra of the operators defined above have the following limits, where $\mathcal{L}_*$ is the linear operator for the planar spiral wave on $L^2\left(\mathbb{R}^2\right)$ with domain $H^2\left(\mathbb{R}^2\right)$ defined in (\ref{eqn:planar_lin_op}), and $\mathcal{L}_{\text{bdy}}$ is the boundary sink defined in (\ref{eqn:boundary_sink_plane}):
\begin{enumerate*}
	\item Bounded disk:  \vspace{-0.5cm}
\begin{align*} \Sigma \left( \mathcal{L}_{*,R}\right) \rightarrow \Sigma_{\text{ext}} \left(\mathcal{L}_*\right) \cup \Sigma_{\text{ext}}\left(\mathcal{L}_{\text{bdy}}\right) \cup \Sigma_{\text{abs}} \left(\mathcal{L}_*\right), \ \text{ as } R \rightarrow \infty \end{align*}

	\item Non-reflecting disk: \vspace{-0.5cm}
\begin{align*} \Sigma \left( \mathcal{L}_{\text{R,nr}}\right) \rightarrow \Sigma_{\text{ext}} \left(\mathcal{L}_*\right)\cup \Sigma_{\text{abs}} \left(\mathcal{L}_*\right), \ \text{ as } R \rightarrow \infty \end{align*}

	\item Boundary sink: \vspace{-0.5cm}
\begin{align*} \Sigma  \left( \mathcal{L}_{\text{bdy,L}}\right) \rightarrow \Sigma_{\text{ext}}\left(\mathcal{L}_{\text{bdy}}\right) \cup \Sigma_{\text{abs}} \left(\mathcal{L}_*\right), \ \text{ as } L \rightarrow \infty.   \end{align*}
\end{enumerate*}

\end{lemma}

\textbf{Numerical Methods:} \\
The patterns and spectra of each operator are computed numerically in Matlab. Patterns are formulated as roots of equations of the form
$\mathcal{F}(U) = 0$ representing the discretized PDE posed on an appropriate domain. Solutions are found using Matlab's built-in root finding algorithm \texttt{fsolve}.
 
Periodic wave trains $U_{\infty}(\xi )$ are found by solving 
\begin{align*}
0 = \kappa^2 D U_{\xi \xi} + \omega U_{\xi} + F(U)
\end{align*}
on the domain $\xi \in [0,2\pi)$ with periodic boundary conditions $U_{\infty}(\xi + 2\pi) = U_{\infty}(\xi)$. Translational symmetry creates a family of solutions, and to select a unique solution and create a square system the phase condition
\begin{align} \label{eqn:phase_condition}
0 = \int_0^{2\pi} \langle U_{\xi}(y), U_{\text{old}}(y) - U(y) \rangle \ d y
\end{align}
and is added to $\mathcal{F}(U)$ where $U_{\text{old}}(\xi)$ is the initial guess for the wave train. One-dimensional periodic domains are discretized using Fourier spectral differentiation matrices with $N_{\xi} = 128$ grid points. Continuous spectra of the wave train are calculated through numerically continuation of the linear dispersion relation $\mathcal{L}^{\text{mv}}_{\infty}(\lambda,\nu)\tilde{V} = 0$ using methods described in \cite{Rademacher:2007uh}, which gives $\Sigma_{\text{FB}}$ of spiral waves via the relation (\ref{eqn:ess_spec_wt_spiral_relation}).
 
On large bounded disks, spiral waves are computed as roots of the equation
\begin{align*}
0 = D \left( \partial_{rr} + \frac{1}{r} \partial_r + \frac{1}{r^2} \partial_{\psi \psi} \right) U + \omega U_{\psi} + F(U)
\end{align*}
with appropriate boundary conditions (homogeneous Neumann or non-reflecting). The spiral angular frequency $\omega_R$ depends on the radius $R$ of the disk and is added as a free parameter in the spiral calculation. Rotational symmetry also creates a family of solutions, and the phase condition for one dimensional waves (\ref{eqn:phase_condition}) is applied at $r = R/2$ to fix the phase of the wave and select a unique solution. 
 
Operators for disk domains $B_R(0)$ and $B_{R}^{\text{nr}}(0)$ are discretized with a fourth-order centered finite difference scheme with $N_r$ = 200 grid points in the radial direction and periodic Fourier spectral methods with $N_{\theta}$ = 100 grid points in the angular coordinate. Grid sizes and discretizations are chosen to ensure numerical accuracy and to capture a sufficient number of spiral bands for convergence of eigenvalues to occur, while maintaining efficient calculations. Radii of $R = 125$ and $R = 5$ are used for the R\"{o}ssler and Karma models, respectively, which results in capturing at least three spiral bands due to the spatial wave numbers. A Neumann compatibility condition is enforced at the origin of the polar grid. As in \cite{Wheeler:2006eu}, two variations of polar grids are used for the spiral and eigenvalue computations. Spiral solutions are solved on a grid of size $N_{\theta} \times N_r$, where the origin contains $N_{\theta}$ grid points. A grid with only one grid point at the origin is used for eigenvalue calculations. Neumann boundary conditions on the outer radius are implemented into the finite difference matrices via the ghost point method \cite[Section 1.4]{Thomas:1995iv}.
%
 
The boundary sink operator on the rectangular domain $\Omega_{\text{bdy}}$ was discretized similarly using fourth-order centered finite differences with $N_s$ grid points in the spatial direction and a Fourier spectral method with $N_t$ grid points in the periodic temporal direction. Following the methods and terminology in \cite{Lloyd:2017jn, Goh:2017sh}, the boundary sink is computed numerically by decomposing the domain into a \enquote{far-field} region in which the asymptotic wave train is translated in time and space and a \enquote{core} region where the Neumann boundary condition has an effect on the wave shape. Note that in this case the core refers to the area near the boundary. The pattern and spatial wave number of the far-field region is fixed to match that of the spiral wave. Smooth cut-off functions of the form $\chi(x) = 1/2\left(1 + \tanh(x - d) \right)$ match the solutions $U_{\text{wt}}(x,\tau)$ in the far-field  and in the core $W(x,\tau)$. The full solution is then given by
\begin{align*}
U_{bdy}(x,\tau) = \left( 1 - \chi(x) \right)U_{wt}(x,\tau) + \chi(x) W(x,\tau).
\end{align*}
Substituting the form of $U_{bdy}$ into (\ref{eqn:bndry_sink}) allows us to calculate $W(x,\tau)$ with Newton's method. To account for numerical inaccuracies, the temporal frequency $\omega$ is set as a free parameter and an integral phase condition is added to match the $U_{wt}$ and $W$ solutions. That is, computing boundary sinks amounts to solving the system
\begin{align*}
&-\omega \partial_{\tau} U_{\text{bdy}} + D \partial_{xx} U_{\text{bdy}}+ F\left( U_{\text{bdy}}\right) = 0,  \ \ (x,\tau) \in (-L,0) \times S^1 \\
&  W_x(0,\tau) = 0, \ \ \tau \in S^1\\
&\int_{-2\pi/\kappa}^0 \int_0^{2\pi} \partial_{x} U_{wt}(x,\tau) W(x,\tau) \ d\tau \ dx = 0
\end{align*}
for $W(x,\tau)$ where the temporal direction is scaled to be $2\pi$-periodic in $\tau$.
 
Solutions in the far-field are obtained by translating asymptotic wave trains in time and space using the angular frequency and spatial wave number from the spiral and imposing the cutoff function $\chi (x)$. Applying $\left(1 - \chi (x)\right)$ to the translation yields an initial condition for $W(x,\tau)$. Domain sizes were selected to fit 6 periods of the wave train, which accurately captured both the Neumann boundary conditions and convergence to the far-field dynamics. Asymptotic wave trains were computed from the one-dimensional problem (\ref{eqn:wt}) using Fourier spectral methods on a periodic grid of $N_t$ points. The translation of wave train to boundary sink resulted in $N_s = 6 N_s$. To take spatial derivatives, the pattern was initially posed on a larger spatial grid of 8 periods with Neumann boundary conditions on each end. When solving for the final pattern, the left two periods were removed to eliminate left-hand side boundary effects and simulate a half-infinite line.

\section{Results}

Spirals on each domain are numerically calculated for the Karma and R\"{o}ssler models, and the influence of the spiral regions is determined by comparing the spectra of the three operators $\mathcal{L}_{*,R}$, $\mathcal{L}_{\text{R,nr}}$, and $\mathcal{L}_{\text{bdy}}$.

\subsection{R\"{o}ssler Model: Line defects are driven by the boundary}

\begin{figure}[H]
\centering
 \includegraphics[width=0.9\linewidth]{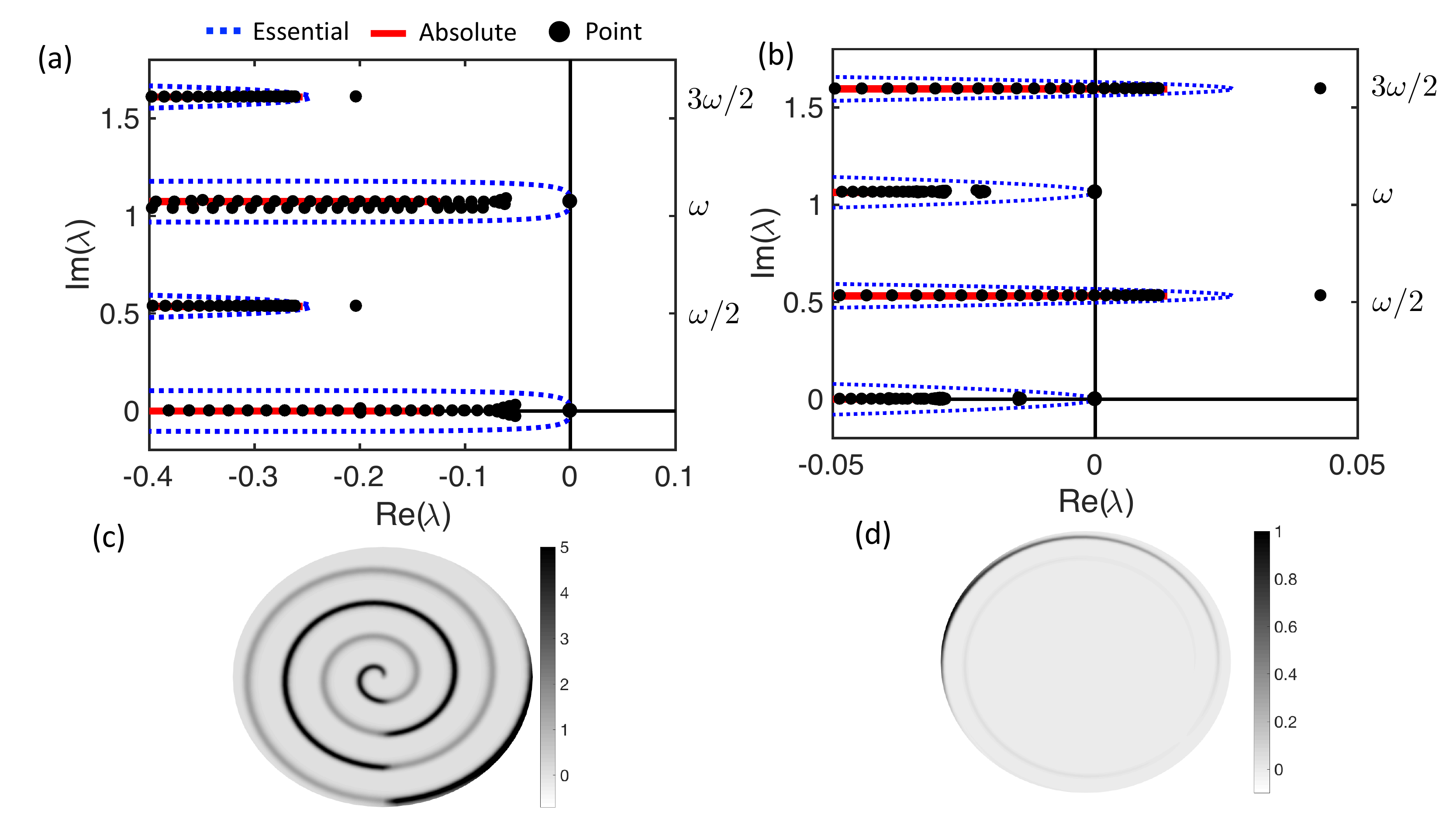}
\caption{ R\"{o}ssler Model: (a) Spectra for $\mathcal{L}_{*,R}$ representing a stable spiral on a disk with parameter $\mu_R= 2$ and radius $R = 125$. Labels on right side of imaginary axis indicate half-multiples of angular frequency. (b) Spectra of unstable spiral, $\mu_R = 3.4$, (c) Spiral on bounded disk of radius $R = 125$ exhibiting a single stationary line defect. Parameter $\mu_R = 3.4$. (d) Unstable point eigenfunction responsible for line defects with $\mu_R = 3.4$. Corresponds to eigenvalue $\lambda = 0.043 + 0.54i = 0.043 + \omega/2 i$.  } \label{fig:fig2}
\end{figure}

\begin{figure}[H]
\centering
 \includegraphics[width=0.85\linewidth]{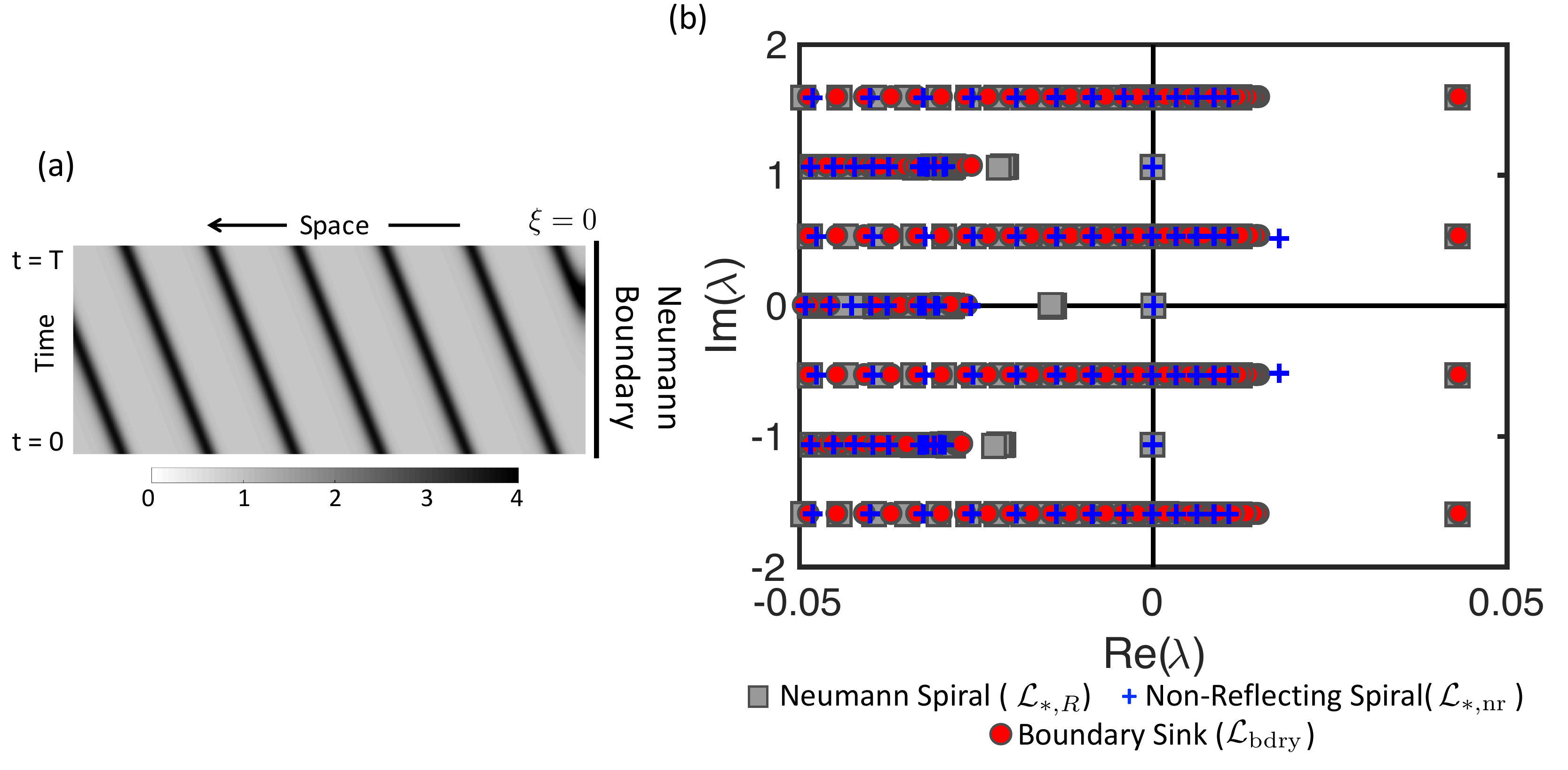}
\caption{ R\"{o}ssler Model: (a) Image of boundary sink. Neumann boundary conditions on the right at $\xi = 0$. Domain is periodic in time (vertical) direction. (b) Spectra of operators $\mathcal{L}_{*,R}$,  $\mathcal{L}_{R,\text{nr}}$, $\mathcal{L}_{\text{bdy}}$.} \label{fig:fig4}
\end{figure}

At the onset of period doubling, point eigenvalues with imaginary parts approximately equal to $\frac{\omega}{2}+\ell\omega$,~$\ell\in\mathbb{Z}$ destabilize, followed by branches of essential and then absolute spectra upon increasing $\mu_R$ further. The unstable eigenfunctions are localized at the boundary (Figure~\ref{fig:fig2}), indicating that line defects are a result of instabilities of the boundary conditions. The spectra of $\mathcal{L}_{*,R},  \mathcal{L}_{\text{R,nr}}$, and $\mathcal{L}_{\text{bdy}}$ are compared in Figure~\ref{fig:fig4}. As expected, all patterns have eigenvalues from the far-field dynamics aligning along the absolute spectrum. However, only domains with boundary conditions, that is the spiral on $B_R(0)$ and boundary sink on $\Omega_{\text{bdy}}$ contain the unstable line defect eigenvalues. Thus, the instability is confirmed to arise from the boundary conditions, and a bounded domain is necessary for the defects to occur.

To further probe for influence of the boundary conditions, we can modify them by changing $\kappa$ in equation (\ref{eqn:asymp_spiral}); $\kappa = 0$ corresponds to homogeneous Neumann conditions and $\kappa = \kappa_*$ is non-reflecting. Therefore, we can start at $\kappa = 0$ with the homogeneous Neumann boundary spiral from $B_R(0)$, numerically continue in $\kappa$ until reaching the spatial wave number of the spiral $\kappa_*$ and track the evolution of an unstable eigenvalue. The spiral is already formulated as a root finding problem, and the eigenvalue problem can be as well by solving $\mathcal{L}_{\text{R,nr}} V - \lambda V = 0.$ Each continuation stage is a 2-step process. First, a new spiral with updated boundary conditions is computed, and second the linearized operator $\mathcal{L}_{\text{R,nr}}$ is modified and an eigenpair $(\lambda,V)$ is computed. Starting the continuation at $\kappa = 0$ allows the unstable eigenfunction from $\mathcal{L}_{*,R}$ to be used in the first continuation step. If the eigenvalue is unchanged with the boundary continuation, then it does not originate from the boundary sink.

\begin{figure}[H]
\centering
 \includegraphics[width=0.9\linewidth]{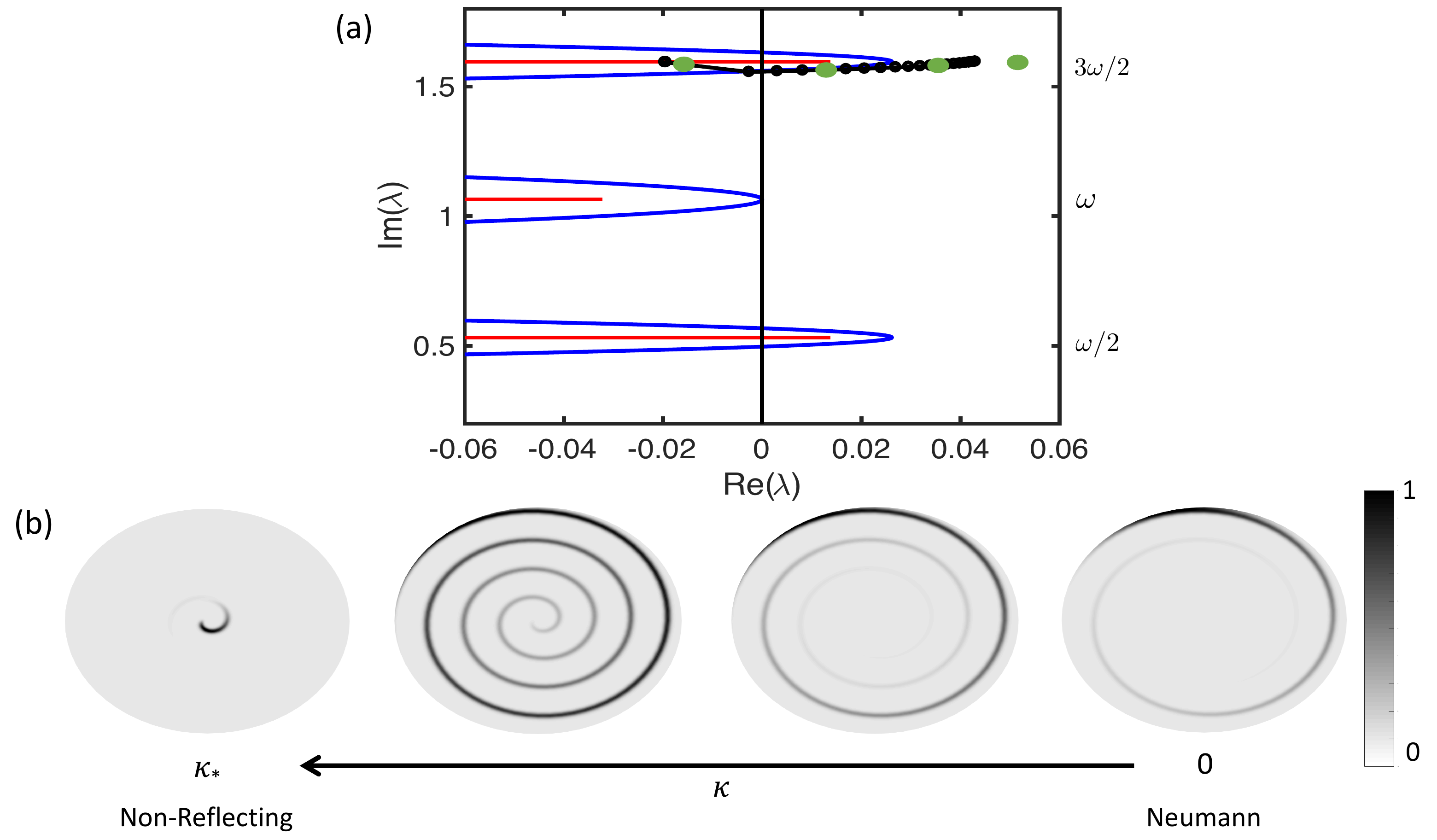}
\caption{ R\"{o}ssler Model: (a) Boundary condition continuation of line defect eigenvalue (black curve), starting with $\kappa = 0$ on the right and ending on absolute spectrum with $\kappa = \kappa_*$.  (b) Eigenfunctions from continuation. Locations of eigenvalues indicated by green circles on eigenvalue continuation in (a). Eigenfunctions are on a disk of radius $R = 125$. } \label{fig:fig5}
\end{figure}

Figure~\ref{fig:fig5} shows the evolution of the point eigenvalue during the $\kappa$-continuation. The eigenvalue changes with $\kappa$, first tracking along the essential spectrum, and then jumping on the absolute spectrum. Increasing $\kappa$ corresponds to a mixed boundary condition and results in different shapes of unstable eigenfunctions, demonstrating that the boundary conditions will change the observed instability. The eigenfunctions in Figure~\ref{fig:fig4}(b) show the transition from localization at the boundary to localization at the core as $\kappa$ is increased from 0 to $\kappa_*$. Similar results are obtained for unstable point eigenvalues at other multiples of $i\frac{\omega}{2} + i \omega\ell$ as these eigenvalues arise due to the asymptotic $e^{i\ell \psi} V(r,\psi)$ symmetry of the eigenfunctions.

\subsection{Karma Model: Alternans are driven by the core}

\begin{figure}[ht]
\centering
 \includegraphics[width=0.9\linewidth]{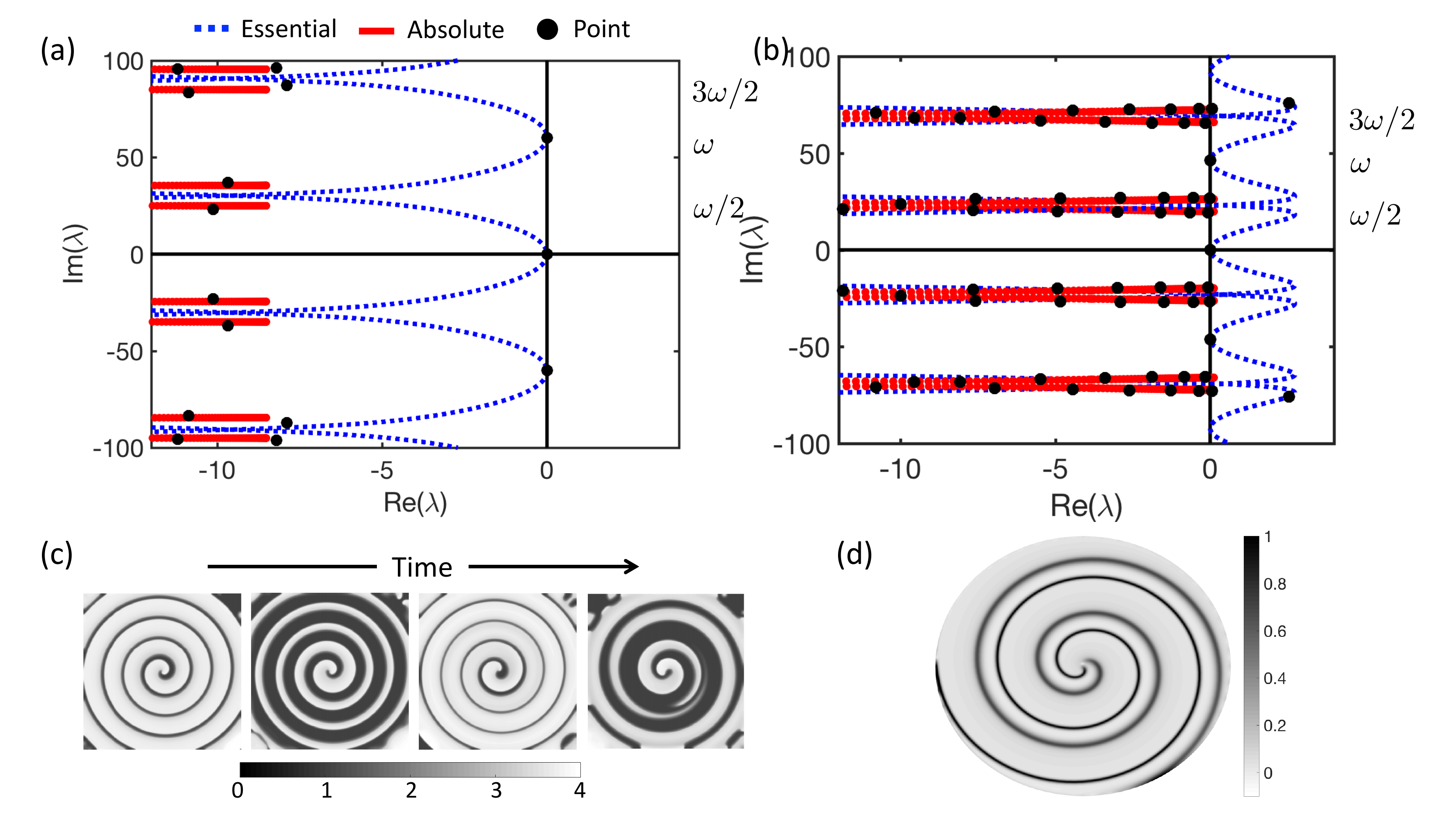}
\caption{ Karma Model: (a) Spectra for stable spiral, $\mu_K= 0.6$. Labels on right side of imaginary axis indicate half-multiples of angular frequency. (b) Spectra of unstable spiral, $\mu_K = 1.4$, (c) Development of alternans in time evolution of spiral on bounded square with homogeneous-Neumann boundary conditions. Parameter $\mu_K = 1.4$. Square of side length 16cm. (d) Unstable point eigenfunction responsible for alternans. Corresponds to eigenvalue $\lambda = 2.6 + 75.9i \approx 2.6 + 3 \omega/2 i$, $\mu_K = 1.4$. Domain radius $R = 5$ with homogeneous-Neumann boundary conditions. } \label{fig:fig3}
\end{figure}

As the bifurcation parameter $\mu_K$ is increased above one in the Karma model, the essential spectrum destabilizes in an Eckhaus instability, followed by a single complex-conjugate pair of eigenvalues with imaginary part near $3\omega/2$. Meandering and alternans appear with the Hopf bifurcation from the point eigenvalues. Shown in Figure~\ref{fig:fig3}d, the unstable point eigenfunction, and hence the form of the instability, has highest magnitude at the boundary of the spiral bands, leading to the observed alternans.  

\begin{figure}[ht]
\centering
 \includegraphics[width=0.9\linewidth]{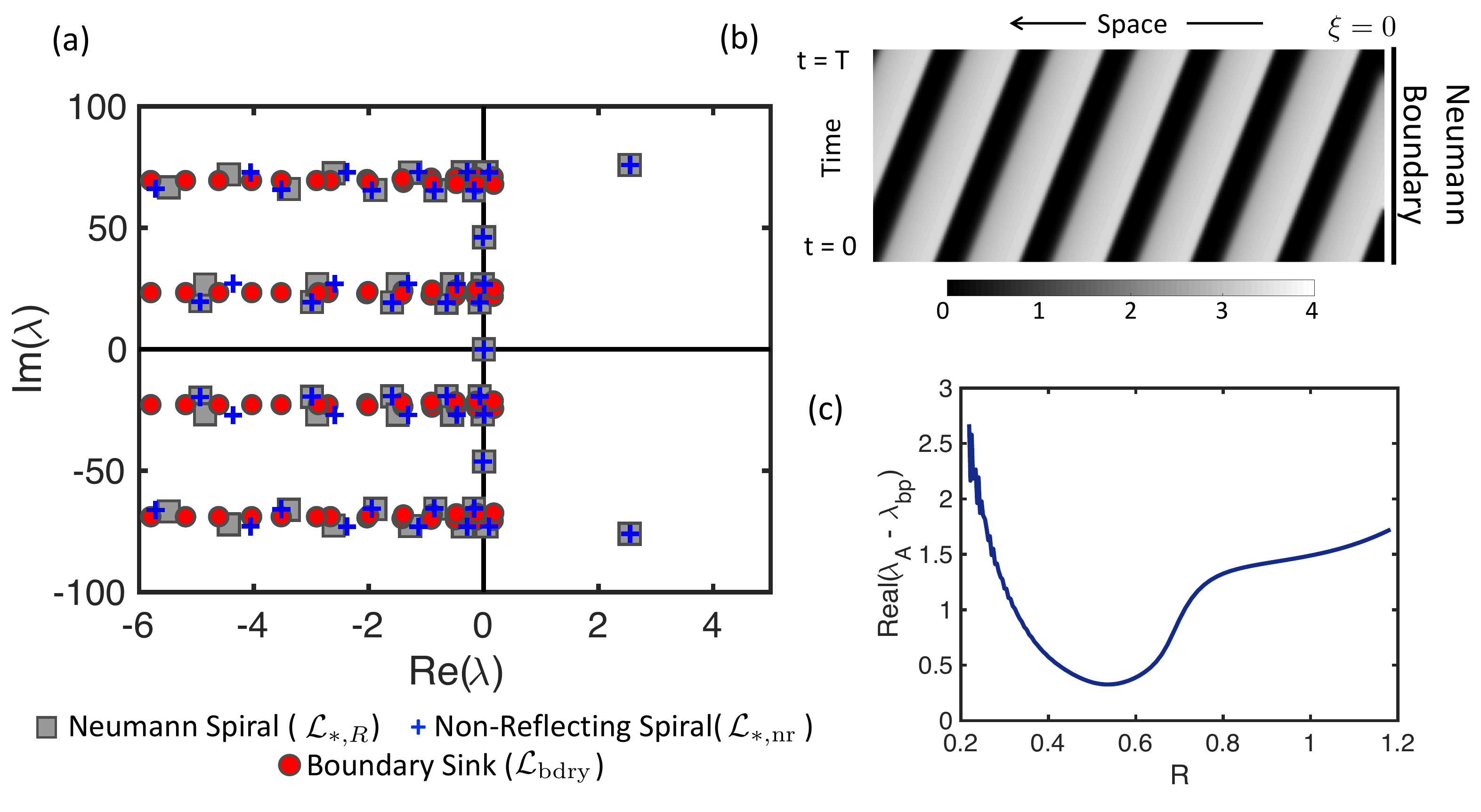}
\caption{  Karma Model: (a) Spectra of operators $\mathcal{L}_{*,R}$,  $\mathcal{L}_{R,\text{nr}}$, $\mathcal{L}_{\text{bdy}}$. (b) Image of boundary sink. Neumann boundary conditions on right at $\xi = 0$. Periodic in time (vertical) direction. (c) Distance between real part of alternans point eigenvalue and branch point of absolute spectrum as a function of parameter $\mu_R$.  } \label{fig:fig6}
\end{figure}

The single pair of eigenvalues suggests they arise from $\Sigma_{\text{ext}}$ and are instabilities of the core. In this case, comparison of the three operators indicates alternans eigenvalues are present in the spiral spectra for $\mathcal{L}_{*,R}$ and $\mathcal{L}_{\text{R,nr}}$, but are absent in the boundary sink $\mathcal{L}_{\text{bdy}}$. Modification of the boundary conditions in $B_{\text{R,nr}}(0)$ by continuing $\kappa$ results in no change to the alternans eigenvalue or eigenfunction. Furthermore, the pair is not emitted from the absolute spectrum; continuation of the absolute spectrum branch point $\lambda_{bp}$ and alternans eigenvalue $\lambda_A$ in parameter $\mu_K$ shows the difference $\text{Re}(\lambda_{bp}) - \text{Re}(\lambda_{A})$ is positive over an appropriate range of parameter values (Figure~\ref{fig:fig6}c). Thus, the point eigenvalues causing the alternans instability stem instead from the unstable pair of eigenvalues originating from $\Sigma_{\text{ext}}$ affiliated with the core.

\subsection{Alternans from interaction of point and essential spectrum}

\begin{figure}[H]
\centering
 \includegraphics[width=0.9\linewidth]{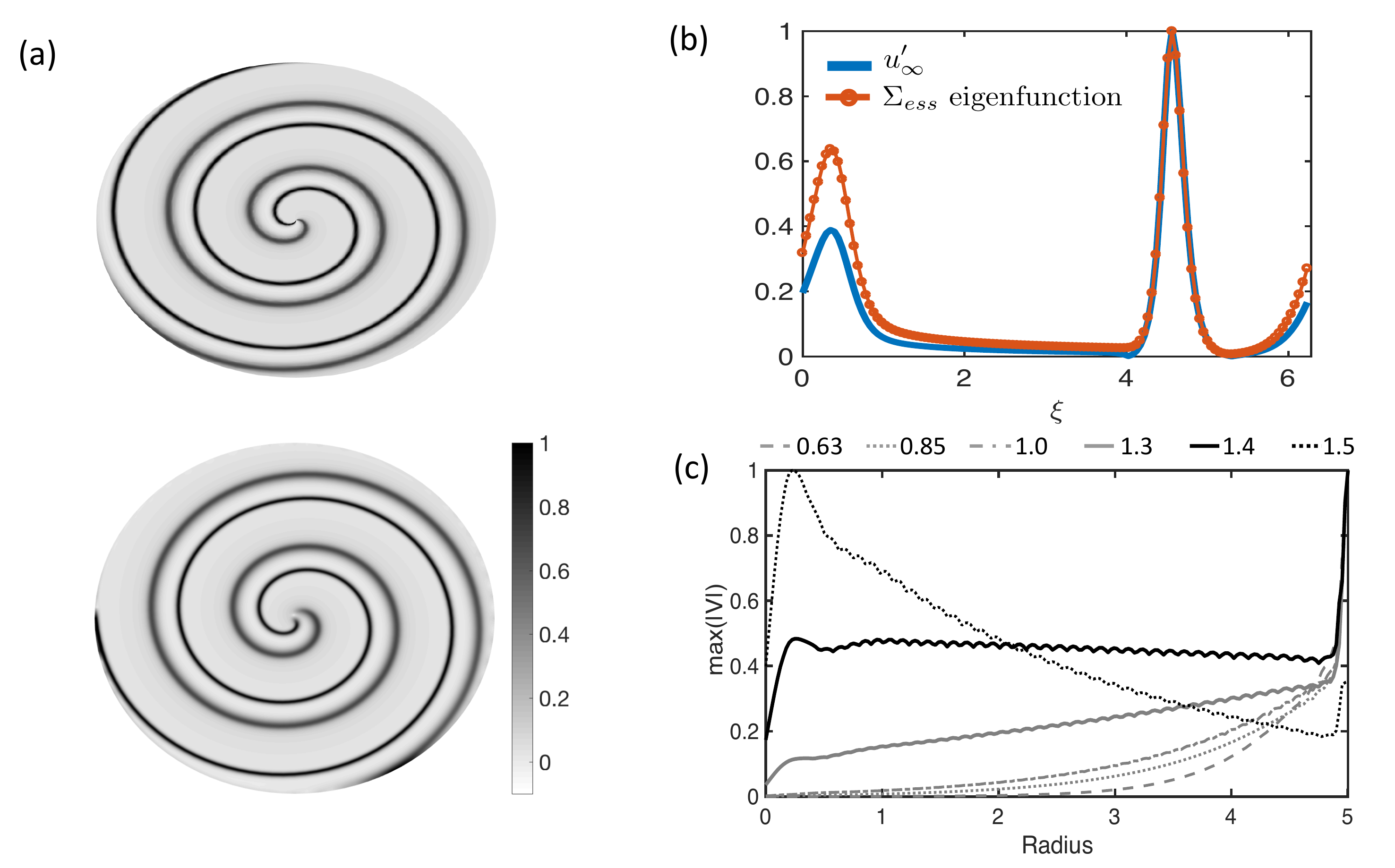}
\caption{  (a) Karma model alternans eigenfunction constructed from essential spectrum data. (b) Comparison between derivate of wave train (blue) and essential spectrum eigenfunction (orange dots), (c) Radial growth of alternans eigenfunctions. Legend indicates value of bifurcation parameter $\mu_K$.} \label{fig:fig7}
\end{figure}

More can be said about the alternans eigenfunction. To leading order, spiral eigenfunctions are of the form (\ref{eqn:efcn_form}), and as 
the unstable alternans point eigenvalue passes through the essential spectrum, the eigenfunction inherits properties of the continuous spectrum  \cite{Sandstede:2001bb}. In particular,  Re $\nu$ is small when the eigenvalue is near the essential spectrum, with Re $\nu$ positive (negative) for the eigenvalue to the left (right) of the continuous spectrum curve. Small Re $\nu$ indicates little radial growth, and Im $\nu = \gamma$ is set by the essential spectrum point. The wave train eigenfunction is determine by the eigenfunction on the essential spectrum, $V_{\text{ess}}$.  Using these observations, the spiral eigenfunction to leading order in the radius is given by
\begin{align}
V(r,\psi) =  e^{i \gamma r} e^{\ell \psi} V_{\text{ess}}(\kappa r + \psi).
\end{align}
A numerically constructed eigenfunction is shown at the top of Figure~\ref{fig:fig7}a. Note that the derivative of the wave train $\partial_{\xi} U_{\infty}(\xi)$ is the eigenfunction on the essential spectrum with eigenvalue $\lambda = i\ell\omega$, $\ell \in \mathbb{Z}$. The alternans eigenfunction crosses $\Sigma_{\text{ess}}$ near one of these points, and therefore $V_{\text{ess}}(\xi)$ is close to $U'_{\infty}(\xi)$ (Figure~\ref{fig:fig7}b). The derivative of the wave train is the highest at the wave fronts and backs and it is this shape that leads to the changing width of the spiral bands and form of alternans. Moreover, the structure of the constructed eigenfunction is in good agreement with the alternans eigenfunction from $\mathcal{L}_{*,R}$ which is reproduced on the bottom of Figure~\ref{fig:fig7}a.

When the alternans eigenvalue is to the left of the essential spectrum, approximately $\mu_K < 1.4$, the overall shape of the eigenfunction is comparable to those shown in Figure~\ref{fig:fig7}a, but there is slight radial growth toward the boundary, corresponding to a spatial eigenvalue, $\nu$, with a small positive real part. Radial growth of the alternans eigenfunction for several values of bifurcation parameter $\mu_K$ is visible in Figure~\ref{fig:fig7}c. As the eigenvalue approaches the essential spectrum, radial growth decreases and is near zero for $\mu_K$ near 1.4, as seen in the solid black curve in Figure~\ref{fig:fig7}c. Alternans appear when this eigenvalue first destabilizes, but the strength of the instability on the core and far-field regions is determined by the proximity of the point eigenvalue to the essential spectrum. 

\section{Discussion}
Instabilities in patterns observed on finite domains may be initiated by unstable eigenvalues stemming from a variety of sources. Determining where unstable eigenvalues originate from yields insight into what creates the accompanying instabilities and what their spatial shape looks like. Here, we present a methodology for unfolding the origin of these eigenvalues and apply it to the specific case of period-doubling instabilities of spiral waves on bounded disks. The technique of comparing the spectra of the three operators can be applied more generally to pattern forming systems on any domain, as long as the patterns of interest can be computed as roots of an appropriate system.
 
Our results predict that line defects in the R\"{o}ssler model will only be seen on bounded domains and that the shape and type of boundary conditions will likely affect the structure of the instability. Furthermore, the interaction of the outer bands of multiple spirals can induce a non-trivial boundary condition between the spirals, and line defects or similar structures may be generated in these situations. Therefore, for instabilities of the boundary sink, an accurate representation of the boundary conditions and domain is a necessary factor when matching models and theory with experiments.

We find that alternans are a product of unstable eigenvalues in the extended point spectrum associated with the spiral core, implying that, as long as the core does not directly interact with the boundary, the shape of the domain and the precise form of boundary conditions are insignificant factors in the spiral stability and formation of alternans. This result has direct impacts for cardiac dynamics in that conclusions for the development of alternans on bounded disks can be extended to irregular and complex geometries such as the heart. Unstable alternans eigenfunctions do exhibit slight radial growth or decay depending on the value of parameter $\mu_K$, which may influence which regions of the spiral are impacted the most by the instability or how prevalent alternans are on a small domain. Our results are consistent with \cite{Marcotte:2015ex,Marcotte:2016hn} which finds that alternans development is most sensitive to perturbations near the core. Both results also suggest that a perturbation must be applied to the core region to destabilize an alternans spiral.
 
Despite the bounded domain, the essential spectrum provides useful information about the source of the alternans instability. Alternans on a ring were previously attributed to destabilization of the continuous spectra \cite{Bar:2004jc}. We find that, on a bounded disk, an unstable point eigenvalue passing through the Eckhaus unstable essential spectrum is fundamental to the alternans structure. Translational symmetry of the wave train ensures existence of the eigenvalue-eigenfunction pair ($\lambda = 0, V = U_{\infty}')$, meaning that the essential spectrum curve that passes through the origin for one parameter set must contain the origin for all parameter values. Therefore, these branches generically destabilize through an Eckhaus instability. Point eigenvalues that interact with these essential spectrum branches acquire an eigenfunction with shape close to $U'_{\infty}$ which impacts the wave fronts and backs leading to the observed form of the alternans instability. The unstable essential spectrum also implies that the associated wave trains on $\mathbb{R}$ undergo an instability as well. 
 
In the R\"{o}ssler model, point eigenvalues do not interact with unstable essential spectra branches and unstable continuous spectra branches do not pass through the origin. Our results suggest that alternans instability will appear if a point eigenvalue crosses essential spectrum that has destabilized through an Eckhaus instability. Future work includes investigation into whether an Eckhaus instability is necessary or sufficient for the formation of alternans. If the continuous spectrum can be attributed to formation of alternans, the 1D computation of the essential spectrum provides a tractable tool for analysis of more realistic and complex models. 
 
There are a number of differences between the spectra of the two models. In the R\"{o}ssler system, countably many discrete eigenvalues destabilize ahead of distinct period-doubling essential and absolute spectrum branches. In contrast, it is a single pair of unstable complex conjugate eigenvalues with imaginary part approximately $3\omega/2$ that leads to alternans. Generically, continuous spectra curves are symmetric around the lines $i\frac{\omega}{2}  + i \omega \ell$, $\ell \in \mathbb{Z}$ which may take the form of smooth curves that intersect with or coincide with symmetry lines, or disjoint branches that do not intersect \cite{Sandstede:2007jw}. The first case is observed in the R\"{o}ssler system, and latter in Karma. Furthermore, in the R\"{o}ssler system, the point $\text{Im}(\lambda) = \frac{\omega}{2}  +  \omega \ell$ on the continuous spectra has spatial eigenvalue $\text{Im}(\nu) =  \kappa/2$, which corresponds with robust period-doubling of the far-field dynamics \cite{Sandstede:2007jw}. On bounded domains discrete eigenvalues limit to absolute spectra curves, meaning an unstable absolute spectra for the R\"{o}ssler system will result in instabilities with temporal frequencies precisely $\omega/2 + \ell \omega$. On the other hand, the disjoint absolute spectrum branches in the Karma model result in an unstable absolute spectrum contributing many frequencies close to, but not specifically period-doubling.
 
In real cardiac systems, alternans lead to spiral break up, providing evidence they originate through a subcritical bifurcation \cite{Frame:1998wi,Rosenbaum:1994bt}. Numerical studies are inconclusive and show both immediate break up and short time alternans persistence \cite{Courtemanche:1993jn,Karma:1993tq,Karma:1994gb}. In \cite{Gottwald:2008hi}, alternans are analytically shown to originate in a subcritical Hopf bifurcation, but this analysis is limited as it relies on a specific normal form for systems near a saddle node of a traveling wave which is not satisfied in all alternans generating models. The debate of a sub- versus super-critical bifurcation can be investigated in  models by determining whether an alternans spiral is stable when considered as a time-periodic three-dimensional structure on a bounded disk.

\paragraph{Acknowledgements.}
Dodson was supported by the NSF through grants DMS-1148284 and 1644760.
Sandstede was partially supported by the NSF through grant DMS-1714429.

\section{Appendix}

The standard form of the Karma model in the laboratory frame for $x \in \mathbb{R}^2$ is
\begin{align*}
E_t &= \gamma \Delta E +  \frac{1}{\tau_E}\left( -E + \left(E^* - n^M \right) \left(1 - \tanh(E - E_h)\right) \frac{E^2}{2} \right)\\
n_t &= \delta \Delta n + \frac{1}{\tau_n} \left( \frac{1}{1 - e^{-Re}} \theta_s \left(E - E_n \right) - n \right),
\end{align*}
where $E = E(x,t)$ represents membrane voltage and $n = n(x,t)$ takes the place of a slower gating variable. In the notation of the main paper, the variables $u$ and $v$ represent $E$ and $n$, respectively. The Heaviside function has been replaced by the smoothed function $\theta_s(u) = (1 + \tanh(su))/2$. Full parameter values are given in Table~\ref{table:params}. The bifurcation parameter (typically called the restitution parameter), $\mu_K = Re$ is increased from 0.6 to 1.4 and controls recovery properties of the excitable media. All other parameters are held fixed in our study.
 
In the laboratory frame, the R\"{o}ssler model is given by
\begin{align*}
u_t &= \delta_1 \Delta u -v - w \\ 
v_t &= \delta_2 \Delta v + u + a v\\ 
w_t & = \delta_3 \Delta w + u w - c w + b .
\end{align*}
The bifurcation parameter $c$ is increased from 2 to 3.4, with line defects appearing as $\mu_R = c$ passes through 3. Parameters are given in Table~\ref{table:params}.

\begin{table}[H]
\begin{center}
\begin{tabular}{ l | l  }
Karma & R\"{o}ssler\\
\hline
$\gamma = 1.1 $ & $\delta_1 = 0.4$\\
$\delta = 0.1$ & $\delta_2 = 0.4$\\
$\tau_E =  0.0025$ & $ \delta_3 = 0.4$\\
$\tau_n = 0.25$ & $ a = 0.2$\\
$E^* = 1.5414$  & $b = 0.2$\\
$M = 4$ & $ \mu_R = c \in [2, 3.4]$\\
$s = 4$ &$\omega \in [1.08, 1.06]$\\
$E_h = 3$ & \\
$E_n = 1 $ & \\
$\mu_K = Re \in [0.6,1.4]$ & \\
$\omega \in [60.02, 46.13]$ & \\
\end{tabular}
\caption{Model parameters: Angular frequency $\omega$ is selected by spirals and given intervals align with bifurcation parameter. } \label{table:params}
\end{center}
\end{table}




\bibliography{period_doubling_instabilities_bib}{}
\bibliographystyle{abbrv}

\end{document}